\RequirePackage{fix-cm}
\documentclass[smallextended]{svjour3}       % onecolumn (second format)
\smartqed  % flush right qed marks, e.g. at end of proof

\usepackage{amsfonts}
\usepackage{amsmath}
\usepackage{amssymb}
\usepackage{graphicx}
\usepackage{multirow}
\usepackage{color}
\usepackage{ulem}
\begin{document}

\title{The $q$-gradient method for global optimization}

\author{Aline C. Soterroni \and
Roberto L. Galski \and
Fernando M. Ramos}

\institute{Aline C. Soterroni \and Roberto L. Galski \and Fernando M. Ramos \\
              Applied Computing Program \and Satellite Control Centre \and Laboratory of Computing and Applied Mathematics \at
              National Institute for Space Research, S\~{a}o Jos\'{e} dos Campos, 12227-010, Brazil \\
              \email{alinecsoterroni@gmailcom, galski@ccs.inpe.br, fernando.ramos@inpe.br}
}

\maketitle

\begin{abstract}
The $q$-gradient is an extension of the classical gradient vector based on the concept of Jackson's derivative. Here we introduce a preliminary version of the $q$-gradient method for unconstrained global optimization. The main idea behind our approach is the use of the negative of the $q$-gradient of the objective function as the search direction. In this sense, the method here proposed is a generalization of the well-known steepest descent method. The use of Jackson's derivative has shown to be an effective mechanism for escaping from local minima. The $q$-gradient method is complemented with strategies to generate the parameter $q$ and to compute the step length in a way that the search process gradually shifts from global in the beginning to almost local search in the end. For testing this new approach, we considered six commonly used test functions and compared our results with three Genetic Algorithms (GAs) considered effective in optimizing multidimensional unimodal and multimodal functions. For the multimodal test functions, the $q$-gradient method outperformed the GAs, reaching the minimum with a better accuracy and with less function evaluations.

\end{abstract}

%%===============================================================================
\section{Introduction}
\label{sec:intro}
%%===============================================================================
Over the last decades the $q$-calculus has been connecting mathematics and physics in applications that span from quantum theory and statistical mechanics, to number theory and combinatorics (see \cite{ernst1} and references therein). Its history dates back to the beginnings of the last century when, based on pioneering works of Euler and Heine, the English reverend Frank Hilton Jackson developed the $q$-calculus in a systematic way \cite{ernst2}. His work gave rise to generalizations of series, functions and special numbers within the context of the $q$-calculus \cite{chaundy}. More important, he reintroduced the concepts of the $q$-derivative \cite{jackson1} (also known as Jackson's derivative) and introduced the $q$-integral \cite{jackson2}.

The $q$-derivative of a function $f(x)$ of one variable is defined as
\begin{equation}
 D_q f(x) = \frac{f(qx)-f(x)}{qx-x},
\end{equation}
where $q$ is a real number different from $1$ and $x$ is different from $0$. In the limit of $q \rightarrow 1$ (or $x \rightarrow 0$), the $q$-derivative reduces to the classical derivative.

Let $f(x)=x^n$, for example. In this case, the classical derivative of $f$ is $nx^{n-1}$ and the $q$-derivative is $[n]x^{n-1}$, where $[n]$ is the $q$-analogue of $n$ given by
\begin{displaymath}
 [n]= \frac{q^n -1}{q-1}.
\end{displaymath}
As $q \rightarrow 1$, $[n]$ tends to $n$. This definition is used to calculate the $q$-binomial and establish a $q$-analogue of Taylor's formula that encompasses many results such as the Euler's identities for $q$-exponential functions, Gauss's $q$-binomial formula, Heine's formula for a $q$-hypergeometric function, among others mathematical results \cite{kac}.

Considering arbitrary functions $f(x)$ and $g(x)$, the $q$-derivative operator satisfy the following properties \cite{kac}:
\begin{enumerate}
 \item [1)] The $q$-derivative is a linear operator for any constants $a$ and $b$
 \begin{displaymath}
 D_q(af(x) + bg(x)) = aD_qf(x) + bD_qg(x).
 \end{displaymath}
 \item [2)] The $q$-derivative of the product of $f(x)$ and $g(x)$ is given by
 \begin{displaymath}
 D_q( f(x)g(x) ) = f(qx) D_qg(x) + g(x)D_qf(x)
\end{displaymath}
that, by symmetry, is equivalent to
 \begin{displaymath}
 D_q( f(x)g(x) ) = f(x) D_qg(x) + g(qx)D_qf(x).
\end{displaymath}
 \item [3)] The $q$-derivative of the quotient of $f(x)$ and $g(x)$ is calculated as
 \begin{displaymath}
 D_q \left( \frac{f(x)}{g(x)} \right)= \frac{g(x) D_qf(x) -f(x)D_qg(x)}{g(x)g(qx)}
 \end{displaymath}
or equivalently
 \begin{displaymath}
 D_q \left( \frac{f(x)}{g(x)} \right) = \frac{g(qx) D_qf(x) -f(qx)D_qg(x)}{g(x)g(qx)}.
 \end{displaymath}
\end{enumerate}
The chain rule for $q$-derivatives does not exist, except for a function of the form $f(u(x))$, where $u(x) = \alpha x^{\beta}$, with $\alpha,\beta$ being constants. More details on the properties of $q$-derivatives can be found in \cite{kac}.

Let a general nonlinear unconstrained optimization problem be defined as
\begin{equation}\label{eq:pmin}
 \min F(\mathbf{x}), \quad  \mathbf{x} = (x_1,\ldots,x_i,\ldots,x_n)
\end{equation}
where $\mathbf{x} \in \Re^n$ is the vector of the independent variables and $F:\Re^n \rightarrow \Re$ is the objective function. The steepest descent method (also known as the gradient descent method) uses information on the gradient of the objective function in seeking the optimum. The search direction is given by the negative of the gradient of $F$. This search strategy is an obvious choice since along this direction the objective function decreases most rapidly.

Requiring only information about first-derivatives, the steepest descent method is attractive because of its limited computational cost and storage requirements \cite{gianni}. However, for multimodal functions, unless one knows in advance where to start from, the search procedure frequently gets stuck in one of the local minima. Consequently, the steepest descent method is not recommended for real-world optimization problems that are usually multimodal. Nevertheless, because of its inherent simplicity, it represents a good starting point for the development of more advanced optimization methods.

Here we propose a generalization of the steepest descent method in which the gradient of the objective function is replaced by its $q$-analogue. Accordingly, the search direction is taken as the negative of the $q$-gradient of $F$. For $q = 1$, the here called $q$-gradient method reduces to the classical steepest descent method. In order to evaluate the performance of the $q$-gradient method we consider three unimodal and three multimodal test functions commonly used as benchmarks. We compare our results with those obtained with the Genetic Algorithms (GAs) G3-PCX developed by Deb et al. \cite{deb}, and the SPC-vSBX and SPC-PNX developed by Ballester and Carter \cite{ballester}, which previous studies have shown to be effective in minimizing multidimensional unimodal and multimodal functions.

The rest of the paper is organized as follows. In Section \ref{sec:qgrad} the $q$-gradient vector is defined. In Section \ref{sec:qgradmethod} the strategies to obtain the parameter $q$ and the step length are described. Section \ref{sec:description} shows the computational experiments and Section \ref{sec:results} discusses the results. Finally, in Section \ref{sec:conclusions} some conclusions and future work are presented.

%%===============================================================================
\section{The $q$-Gradient}
\label{sec:qgrad}
%%===============================================================================
Given a differentiable function of $n$ variables $F(\mathbf{x})$, the gradient of $F$ is the vector of the $n$ first-order partial derivatives of $F$. Similarly, the $q$-gradient is the vector of the $n$ first-order partial $q$-derivatives of $F$. Thus, let the parameter $q$ be a vector $\mathbf{q} = (q_1,\ldots,q_i,\ldots,q_n)$, where $q_i \neq 1 \ \forall i$, the first-order partial $q$-derivative with respect to the variable $x_i$ is given by
\begin{equation}\small
D_{q_i,x_i} F(\mathbf{x}) = \displaystyle \frac{F(x_1,\ldots,q_ix_i,\ldots,x_n)- F(x_1,\ldots, x_i,\ldots,x_n)} {q_ix_i - x_i}
\end{equation}
with
\begin{equation}
 \left. D_{q_i,x_i} F(\mathbf{x}) \right|_{x_i=0} = \frac{\partial F(\mathbf{x})}{\partial x_i}
\end{equation}
and
\begin{equation}
 \left. D_{q_i,x_i} F(\mathbf{x}) \right|_{q_i=1} = \frac{\partial F(\mathbf{x})}{\partial x_i}.
\end{equation}
This framework can be extended to define the $q$-gradient of a function of $n$ variables as
\begin{equation}\label{eq:qgrad}
\nabla_{\mathbf{q}} F(\mathbf{x}) = \left[ D_{q_1,x_1} F(\mathbf{x}) \ \ldots \ D_{q_i,x_i} F(\mathbf{x}) \ \ldots \ D_{q_n,x_n} F(\mathbf{x})\right]
\end{equation}
with the classical gradient being recovered in the limit of $q_i \rightarrow 1$, for all $i=1,\ldots,n$. 

\begin{figure}[ht!]
  \centering
  \includegraphics[width=3.5in]{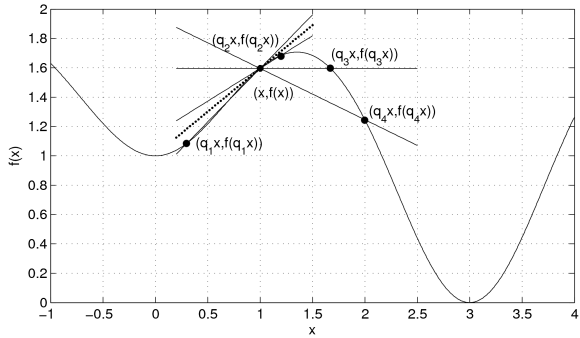}
  \caption{Geometric interpretation of the classical derivative (dotted line) and the $q$-derivative of $f(x)$ at $x = 1.0$ and different values of the parameter $q$.}
  \label{fig:geometrico}
\end{figure}

The Fig. \ref{fig:geometrico} illustrates the geometric interpretation of the classical gradient and the $q$-gradient of a function of one variable $f(x) = 2 -( e^{-x^2} + 2 e^{-(x-3)^2})$. In this case, the gradient is simply the slope of the tangent line at $x$. Similarly, the $q$-gradient is the slope of the secant line passing through the points $(x,f(x))$ and $(qx,f(qx))$. If the slope of the secant line is positive (negative), the $q$-gradient points to the right (left) direction. 

Since the slope of the tangent line (dotted line) at $x=1$ is positive, the steepest descent method at this point will move necessarily to the left and, thus, will be trapped by the local minimum at $x=0$. The slope of the secant line passing through the points $(x,f(x))$ and $(qx,f(qx))$ can be positive or negative depending on the value of the parameter $q$. For instance, if $q=2$ the $q$-derivative is negative at $x=1$ (see the secant line passing through $(x,f(x))$ and $(q_4x,f(q_4x))$ in Fig. \ref{fig:geometrico}), which potentially allows a minimization strategy based on the value of the $q$-gradient to take a leap to the right, towards the global minimum of $f$. Note that there is a value $1.5 < q_3 < 2$ for which $x=1$ is a stationary point of the $q$-gradient ($\nabla_qF(x)=0$) but that does not correspond to any minimum or maximum of $f$. The stationary points of the $q$-gradient are avoided in the method by generating a new the parameter $q$. Finally, for $0<q_1<0.5$ or $1<q_2<1.5$ the slope of the secant line is positive and the $q$-gradient method will move to the left as the steepest descent method.

This simple example above shows that the use of the $q$-gradient offers a new mechanism to escape from local minima. Moreover, the transition from global to local search might be controlled by the parameter $q$, provided a suitable strategy for generating $q$-values is incorporated into the minimization algorithm.

%%===============================================================================
\section{$q$-Gradient Method Description}
\label{sec:qgradmethod}
%%===============================================================================
A general optimization strategy is to consider an iterative procedure that, starting from $\mathbf{x}^0$, generates a sequence $\{\mathbf{x}^k\}$ given by \cite{vanderplaats}
\begin{equation}\label{eq:iterprocess}
 \mathbf{x}^{k+1} = \mathbf{x}^{k} + \alpha^{k} \mathbf{d}^{k}
\end{equation}
where $\mathbf{d}^{k}$ is the search direction and $\alpha^{k}$ is the step length or the distance moved along $\mathbf{d}^{k}$ in the iteration $k$.

Optimization methods can be characterized according to the direction and step length used in Eq. (\ref{eq:iterprocess}). The steepest descent method sets $\mathbf{d}^{k} = -\nabla F(\mathbf{x}^{k})$ and the step length $\alpha^{k}$ is usually determined by a line-search technique that minimizes the objective function along the direction $\mathbf{d}^{k}$. In the $q$-gradient method, as here proposed, the search direction is the negative of the $q$-gradient of the objective function $- \nabla_{q} F(\mathbf{x})$. Thus the optimization procedure defined by Eq. (\ref{eq:iterprocess}) becomes
\begin{equation}\label{eq:search}
 \mathbf{x}^{k+1} = \mathbf{x}^{k} - \alpha^{k} \nabla_{q} F(\mathbf{x}^{k}).
\end{equation}
Key to the performance of the $q$-gradient method, the strategies used to specify the parameter $q$ and the step length $\alpha$ are described below.

%%===============================================================================
\subsection{Parameter $\mathbf{q}$}
\label{sub:parameter}
%%===============================================================================
Considering a function of $n$ variables $F(\mathbf{x})$, a set of $n$ different parameters $q_i \in \Re-\{1\}$ ($i=1,\ldots,n$) are needed to compute the $q$-gradient vector of $F$. The overall strategy adopted here is to draw the values of $q_i$ (or some variable related to them) from some suitable probability density function (pdf), and with a standard deviation that decreases as the iterative search proceeds. In this sense, the role of the standard deviation here is reminiscent of the one played by the temperature in a simulated annealing (SA) algorithm, that is, to make the algorithm go from a very random (at the beginning) to a very deterministic search (at the end).

In the current implementation we opted to first draw the values of $q_i^kx_i^k$ from a Gaussian pdf given by
\begin{equation}
 f(x) = \frac{1}{\sigma \sqrt(2\pi)} \exp\left[ - \frac{(x-\mu)^2}{2\sigma^2} \right],
\end{equation}
with $\mu = x_i^k$ and $\sigma = \sigma^k$; then, we computed the values of $q_i^k$.

Starting from $\sigma^0$, the standard deviation of the pdf is decreased by the following ``cooling'' schedule, $\sigma^{k+1} = \beta \cdot \sigma^{k}$, where $0<\beta <1$ is the reduction factor. As $\sigma^{k}$ approaches zero, the values of $q_i^k$ tend to unity, the algorithm reduces to the steepest descent method, and the search process becomes essentially local. As in a SA algorithm, the performance of the minimization algorithm depends crucially on the choice of parameters $\sigma^0$ and $\beta$. A too rapid decrease of $\sigma^{k}$, for example, may cause the algorithm to be trapped in a local minimum.

%%===============================================================================
\subsection{Step Length}
\label{suc:step}
%%===============================================================================
The calculation of the step length $\alpha$ is a tradeoff. On the one hand, $\alpha$ should give a considerable reduction of the objective function. On the other hand, its calculation should not take too many evaluations of $F$ \cite{nocedal}. Steepest descent algorithms generally use line-search techniques to determine the step length $\alpha^k$ along the steepest descent direction $\mathbf{d}^k=-\nabla F(\mathbf{x}^k)$ at the iteration $k$. A first version of our algorithm \cite{soterroni} applied the golden section for step length determination. However, traditional line-search algorithms, like the golden section, ensure that the condition $F(\mathbf{x}^{k+1}) < F(\mathbf{x}^k)$ is always satisfied, what obviously is a poor strategy when dealing with multimodal minimization problems. In addition, depending on the value of $q$, the negative of the $q$-gradient may not point to the local descent direction.

One way to circumvent these problems is to use a diminishing step length $\alpha^k$, i.e., the initial step length $\alpha^0$ is reduced by $\alpha^{k+1} = \beta \cdot \alpha^{k}$, where, for the sake of simplicity, $\beta$ is the same reduction factor used to compute $\sigma^{k}$. As the step length decreases (and the values of $q_i^k$, in parallel, tend to unity), a smooth transition to an increasingly local search process occurs. \\

%%%===============================================================================
%\subsection{Algorithm for the $q$-Gradient Method}
%\label{sub:algorithm}
%%%===============================================================================
%Based on the definitions presented in the previous sections, the $q$-gradient method for continuous global optimization problems is described as follows. \\
%
%\textsc{Algorithm 1} ($q$-Gradient Method)
%\begin{itemize}
% \item [] Given $\mathbf{x}^0$ (initial point), $\sigma^0>0$, $\alpha^0>0$ and $0< \beta < 1$:
% \item [] 1) Set $k= 0$
% \item [] 2) Set $\mathbf{x}_{best} = \mathbf{x}^k$
% \item [] 3) While the stopping criteria are not reached, do:
%  \begin{itemize}
% \item [] a) Generate $\mathbf{q}^k\mathbf{x}^k$ by a Gaussian distribution with $\sigma^k$ and $\mu^k=\mathbf{x}^k$
% \item [] b) Calculate the $q$-gradient $\nabla_q F(\mathbf{x}^k)$
% \item [] c) Set $\mathbf{d}^{k} = - \nabla_q F(\mathbf{x}^k) / \| \nabla_q F(\mathbf{x}^k) \|$
% \item [] d) Set $\mathbf{x}^{k+1} = \mathbf{x}^k + \alpha^k \cdot \mathbf{d}^{k}$
% \item [] e) If $F(\mathbf{x}^{k+1}) < F(\mathbf{x}_{best})$ set $\mathbf{x}_{best} = \mathbf{x}^{k+1}$
% \item [] f) Set $\sigma^{k+1} = \beta \cdot \sigma^{k}$ and $\alpha^{k+1} = \beta \cdot \alpha^{k}$
% \item [] g) Set $k = k + 1$
%\end{itemize}
% \item [] 4) Return $\mathbf{x}_{best}$.
%\end{itemize}

The main idea behind the $q$-gradient method is to use the negative of the $q$-gradient of $F$, instead of the negative of the classical gradient of $F$, as the search direction. Strategies for generating the parameter $\mathbf{q}^k$ and the step length $\alpha^k$, at each iteration, complement this very simple algorithm. Note that there are three free parameters to be specified, namely, $\sigma^0$, $\alpha^0$ and $\beta$. The initial standard deviation $\sigma^0$ determines how stochastic is the search. For multimodal functions, it must be high enough to allow the method to properly sample the search space. The reduction factor $\beta$ controls the speed of the transition from stochastic to deterministic search. A $\beta$ close to $1$ reduces the risk of being trapped in a local minimum. The last free parameter, the initial step length $\alpha^0$, depends heavily on the topology of the search space and, thus, requires some empirical exploration. In the end, as with the choice of the cooling schedule in a SA algorithm \cite{locatelli}, an appropriate specification of the three free parameters is strictly dependent on the objective function. Although a bad choice may lead to some deterioration in its performance, the $q$-gradient method has shown to be sufficiently robust to still be capable of reaching the global minimum. 

%The algorithm stops when the appropriate stopping criterium is attained. In real-world applications (i.e.,  in problems for which the global minimum is not known), it can be the maximum number of function evaluations, or the value of the local gradient $||\nabla F(\mathbf{x}^k)|| < \epsilon$, since the $q$-gradient method converges to the steepest descent method at the end of the search. The algorithm returns the $\mathbf{x}_{best}$ as the minimum value of the objective function $F$ obtained during the iterative procedure, i.e., $F(\mathbf{x}_{best}) \leq F(\mathbf{x}^k)$, $\forall k$.

%%===============================================================================
\section{Computational Experiments}
\label{sec:description}
%%===============================================================================
The performance of the $q$-gradient method was evaluated over six $20$-variable test functions ($n=20$) commonly employed in the literature. We use the same experimental setup and stopping criteria as described in \cite{deb} and \cite{ballester} in order to allow a direct comparison with their results. The stopping criteria are: maximum of $10^6$ function evaluations or $F(\mathbf{x})< 10^{-20}$. 

As in \cite{deb} and \cite{ballester}, we set the three free parameters for the $q$-gradient method after preliminary exploratory runs. The results presented here are for those which yielded the best performance. The benchmark consists of the following analytical functions:

\begin{enumerate}
 \item [1)] Ellipsoidal function ($F_{elp}$)
\begin{equation}\label{eq:ellipsoidal}
 F_{elp}= \sum_{i=1}^{n} i  \mathbf{x}_i^2.
 \end{equation}
Although the Ellipsoidal function is convex and unimodal, it is an example of a poorly scaled function.
 \item []
 \item [2)] Schwefel's function ($F_{sch}$)
\begin{equation}\label{eq:schwefel}
  F_{sch}= \sum_{i=1}^{n} \left( \sum_{j=1}^{i} \mathbf{x}_j \right)^2.
 \end{equation}
The Schwefel's function is an extension of the Ellipsoidal function and it is also a unimodal and poorly scaled function.
 \item []
 \item [3)] Generalized Rosenbrock's function ($F_{ros}$)
\begin{equation}\label{eq:rosenbrock}
  F_{ros}= \sum_{i=1}^{n-1} [ 100 \cdot (\mathbf{x}_i^2- \mathbf{x}_{i+1})^2 +(1-\mathbf{x}_i)^2].
 \end{equation}
Although the Rosenbrock's function is a well-known unimodal function for $n=2$, numerical experiments have shown that for $ 4 \leq n \leq 30$ the function has two minima, the global one at $\mathbf{x}=\mathbf{1}$ and a local minimum that changes with the dimensionality $n$ \cite{shang}. The Rosenbrock's function is considered a test case for premature convergence once the global minimum lays inside a long, narrow, and parabolic shaped flat valley.
 \item []
 \item [4)] Ackley's function ($F_{ackl}$)
\setlength{\arraycolsep}{0.0em}
\begin{equation}
F_{ackl}=20 + e -20 \exp{\left(-0.2 \sqrt{\frac{1}{n} \sum_{i=1}^{n}\mathbf{x}_i^2}\right)} - \exp{\left(\frac{1}{n} \sum_{i=1}^{n} \cos(2 \pi \mathbf{x}_i )\right)}.
\end{equation}
The Ackley's function is highly multimodal and the basin of the local minima increase in size as one moves away from the global minimum \cite{ballester}.
 \item []
 \item [5)] Rastrigin's function ($F_{rtg}$)
\begin{equation}\label{eq:rastrigin}
F_{rtg}= 10 n + \sum_{i=1}^{n} ( \mathbf{x}_i^2-10\cos(2\pi\mathbf{x}_i) ).
\end{equation}
The Rastrigin's function has a parabolic landscape away from the global minimum, but as we move towards the global minimum, the size of the basins increase. The function is highly multimodal and its characteristics are known to be difficult for many optimization algorithms to achieve the global minimum \cite{ballester}.
 \item []
 \item [6)] Rotated Rastrigin's function ($F_{rrtg}$).
\setlength{\arraycolsep}{0.0em}
\begin{eqnarray}
F_{rrtg}&{}={}& 10 n + \sum_{i=1}^{n} ( \mathbf{y}_i^2 -10\cos(2\pi\mathbf{y}_i) ), \quad \mathbf{y} = A\cdot \mathbf{x} \\
&&{}\: A_{i,i} = 4/5 \nonumber \nonumber \\
&&{}\: A_{i,i+1} = 3/5 \ (i \ \mbox{odd}) \nonumber \\
&&{}\: A_{i,i-1} = -3/5 \ (i \ \mbox{even}) \nonumber  \\
&&{}\: A_{i,j} = 0 \ (\mbox{otherwise}) \nonumber
\end{eqnarray}
The rotated Rastrigin's is a highly multimodal function without local minima arranged along the axis \cite{ballester}.
 \item []
\end{enumerate}
For all these functions the global minimum is $F(\mathbf{x}^{*}) = 0$ at $\mathbf{x}^{*} = \mathbf{0}$, except the Generalized Rosenbrock's function where $\mathbf{x}^{*} = \mathbf{1}$. The initial point set $\mathbf{x}^0$ for each function is generated by a uniform distribution within $[-10,-5]$, as used in \cite{deb} and \cite{ballester}.

%%===============================================================================
\section{Results}
\label{sec:results}
%%===============================================================================

Extensive comparisons between the GAs G3-PCX (results obtained from \cite{deb} for Ellipsoidal, Schwefel, Rosenbrock and Rastrigin functions; and from \cite{ballester} for Ackley and Rotated Rastrigin), SPC-vSBX and SPC-PNX (results obtained from \cite{ballester}) and the $q$-gradient method are presented in Tables \ref{tab:unimodal} and \ref{tab:multimodal}. As in \cite{deb} and \cite{ballester}, the ``Best'', ``Median'' and ``Worst'' columns refer to the number of function evaluations required to reach the accuracy $10^{-20}$. When this condition is not achieved, the best value found so far for the test function after $10^6$ evaluations is given in column ``$F(\mathbf{x}_{best})$''. The column ``Success'' refers to how many runs reached the target accuracy, for unimodal functions, or ended up within the global minimum basin, for multimodal ones. The best performances are highlighted in bold in each table. The corresponding values of the best parameters $\sigma^0$, $\alpha^0$ and $\beta$ used in each test function are given in Table \ref{tab:parameters}.

\begin{table}[ht!]
 \renewcommand{\arraystretch}{1.3}
 \caption{Parameters used by the $q$-gradient method over the test functions.}
 \label{tab:parameters}
 \begin{tabular}{llll}
 \hline\noalign{\smallskip}
 Functions & $\mathbf{\sigma^0} \quad$ & $\mathbf{\alpha^0} \quad$ & $\mathbf{\beta} \quad$ \\
 \noalign{\smallskip}\hline\noalign{\smallskip}
 Ellipsoidal & $0.4$ & $38$  & $0.86$ \\
 Schwefel    & $0.1$ & $1$ & $0.997$ \\
 Rosenbrock  & $0.1$ & $0.1$ & $0.9995$ \\
 \noalign{\smallskip}\hline
 Ackley      & $20$  & $12$  & $0.90$ \\
 Rastrigin   & $21$  & $0.3$ & $0.9995$ \\
 Rotated Rastrigin & $30$  & $0.5$ & $0.999$ \\
 \noalign{\smallskip}\hline
 \end{tabular}
\end{table}

\begin{table}[ht!]
 \renewcommand{\arraystretch}{1.3}
 \caption{Performance comparison between G3-PCX \cite{deb}, SPC-vSBX \cite{ballester}, SPC-PNX \cite{ballester} and $q$-gradient method over the unimodal test functions in terms of the best, median and worst number of function evaluations required to reach the accuracy $10^{-20}$.}
 \label{tab:unimodal}
 \begin{tabular}{lllllll}
 \hline\noalign{\smallskip}
 Function & Method & Best & Median & Worst & $F(\mathbf{x}_{best})$ & Success \\
 \noalign{\smallskip}\hline
 \multirow{4}{*}{Ellipsoidal} & \textbf{G3-PCX} &  $\mathbf{5,826}$  & $\mathbf{6,800}$  & $\mathbf{7,728}$  & $\mathbf{10^{-20}}$ & $\mathbf{10/10}$ \\
 & SPC-vSBX     &  $49,084$ & $50,952$ & $57,479$ & $10^{-20}$ & $10/10$ \\
 & SPC-PNX      &  $36,360$ & $39,360$ & $40,905$ & $10^{-20}$ & $10/10$ \\
& \textbf{$q$-Gradient}  &  $\mathbf{5,905}$ & $\mathbf{7,053}$ & $\mathbf{7,381}$ & $\mathbf{10^{-20}}$ & $\mathbf{50/50}$ \\
 \noalign{\smallskip}\hline
 \multirow{4}{*}{Schwefel} & \textbf{G3-PCX} &  $\mathbf{13,988}$  & $\mathbf{15,602}$  & $\mathbf{17,188}$  &  $\mathbf{10^{-20}}$ & $\mathbf{10/10}$ \\
 & SPC-vSBX     &  $260,442$ & $294,231$ & $334,743$ &  $10^{-20}$ & $10/10$ \\
 & SPC-PNX      &  $236,342$ & $283,321$ & $299,301$ &  $10^{-20}$ & $10/10$ \\
 & $q$-Gradient &  $289,174$  & $296,103$ & $299,178$ &  $10^{-20}$ & $50/50$ \\
 \noalign{\smallskip}\hline
 \multirow{4}{*}{Rosenbrock} & \textbf{G3-PCX} &  $\mathbf{16,508}$  & $\mathbf{21,452}$  & $\mathbf{25,520}$  &  $\mathbf{10^{-20}}$ & $\mathbf{36/50}$ \\
 & SPC-vSBX     &  $10^6$    &     -     &     -     &  $10^{-4}$  & $48/50$ \\
 & SPC-PNX      &  $10^6$    &     -     &     -     &  $10^{-10}$ & $38/50$ \\
 & $q$-Gradient &  $10^6$    &     -     &     -     &  $10^{-10}$ & $50/50$ \\
 \noalign{\smallskip}\hline
 \end{tabular}
\end{table}

\begin{table}[ht!]
 \renewcommand{\arraystretch}{1.3}
 \caption{Performance comparison between G3-PCX \cite{deb,ballester}, SPC-vSBX \cite{ballester}, SPC-PNX \cite{ballester} and $q$-gradient method over the multimodal test functions in terms of the best, median and worst number of function evaluations required to reach the accuracy $10^{-20}$.}
 \label{tab:multimodal}
 \begin{tabular}{lllllll}
 \hline\noalign{\smallskip}
 Function & Method & Best & Median & Worst & $F(\mathbf{x}_{best})$ & Success \\
 \noalign{\smallskip}\hline
 \multirow{4}{*}{Ackley} &  G3-PCX       &  $10^6$   &     -    &     -    & $3,959$    & $0$ \\
  & SPC-vSBX     &  $57,463$ & $63,899$ & $65,902$ & $10^{-10}$ & 10/10 \\
  & SPC-PNX      &  $45,736$ & $48,095$ & $49,392$ & $10^{-10}$ & 10/10 \\
  & \textbf{$q$-Gradient} & $\mathbf{11,850}$ & $\mathbf{12,465}$ & $\mathbf{13,039}$ & $\mathbf{10^{-15}}$ & $\mathbf{50/50}$ \\
 \noalign{\smallskip}\hline
 \multirow{4}{*}{Rastrigin} & G3-PCX       &  $10^6$    &     -     &     -     & $15,936$   & $0$ \\
 &  SPC-vSBX &  $260,685$ & $306,819$ & $418,482$ & $10^{-20}$ & 6/10 \\
 &  SPC-PNX      &  $10^6$    &     -     &     -     & $4.975$    & 0 \\
 &  \textbf{$q$-Gradient} &  $\mathbf{676,050}$  & $\mathbf{692,450}$  & $\mathbf{705,037}$ &  $\mathbf{10^{-20}}$ & $\mathbf{48/50}$ \\
 \noalign{\smallskip}\hline
 &  G3-PCX       &  $10^6$   &     -     &     -    & $309.429$  & $0$ \\
 Rotated   &  SPC-vSBX       &  $10^6$   &     -    &     -      & $8.955$  & $0$ \\
 Rastrigin &  SPC-PNX        &  $10^6$   &     -    &     -      & $3.980$  & $0$ \\
 &  \textbf{$q$-Gradient}    &  $\mathbf{541,857}$  &  $\mathbf{545,957}$ &  $\mathbf{549,114}$   & $\mathbf{10^{-20}}$  & $\mathbf{20/50}$ \\
 \noalign{\smallskip}\hline
 \end{tabular}
\end{table}

In Table \ref{tab:unimodal}, for the Ellipsoidal function, the $q$-gradient method achieved the required accuracy $10^{-20}$ for all 50 runs, with an overall performance similar to the one displayed by the G3-PCX, the best algorithm among the GAs. As for the Schwefel's function, the $q$-gradient method again attained the required accuracy for all runs but was outperformed by the G3-PCX in terms of the number of function evaluations. Finally, for the Rosenbrock's function, the $q$-gradient was beaten by the G3-PCX (the only to achieve the required accuracy) but performed better then the two other GAs. The overall evaluation of the $q$-gradient method performance in these numerical experiments with unimodal (or quasi-unimodal) test functions indicates that it reaches the required accuracy (or the minimum global basin) in $100\%$ of the runs, but it is not faster than the G3-PCX. This picture improves a lot when it comes to tackle the multimodal Ackley's and Rastringin's functions. 

In Table \ref{tab:multimodal}, due to limited computing precision the required accuracy for the Ackley's function was set to $10^{-10}$ for the GAs and $10^{-15}$ in our simulations\footnote{Numerical experiments have shown that Ackley's function evaluated at $\mathbf{x}=0$ with double precision is equal to $-0.4440892098500626$E$-015$ and not zero.}. The $q$-gradient method was here clearly better than the GAs, reaching the required accuracy in more runs or in less functions evaluations. For the Rastrigin's function, the G3-PCX and the SPC-PNX were unable to attain the global minimum basin. The other two algorithms reached the required accuracy $10^{-20}$, but the $q$-gradient method was the only to do it in $96\%$ of the runs (48 over 50). Finally, in the case of the rotated Rastrigin's function, the $q$-gradient was the only algorithm to reach the minimum, attaining the required accuracy in $20$ out of $50$ independent runs. Summarizing the results with multimodal functions, we may say that the $q$-gradient method outperformed the GAs in all the three test cases considered, reaching the minimum with less function evaluations or in more independent runs.

%%===============================================================================
\section{Conclusions and Future Work}
\label{sec:conclusions}
%%===============================================================================
The main idea behind the $q$-gradient method is the use of the negative of the $q$-gradient of the objective function --- a generalization of the classical gradient based on the Jackson's derivative --- as the search direction. The use of Jackson's derivative provides an effective mechanism for escaping from local minima. The method has strategies for generating the parameter $q$ and the step length that makes the search process gradually shifts from global in the beginning to almost local search in the end.

For testing this new approach, we considered six commonly used 20-variable test functions. These functions display features of real-world optimization problems (multimodality, for example) and are notoriously difficult for optimization algorithms to handle. We compared the $q$-gradient method with GAs developed by Deb et al. \cite{deb}, and Ballester and Carter \cite{ballester} with promising results. Overall, the $q$-gradient method clearly beat the competition in the hardest test cases, those dealing with the multimodal functions.

It comes without suprise the (relatively) poor results of the $q$-gradient method with the Rosenbrock's function, a unimodal test function specially difficult to be solved by the steepest descent method. This result highlights the need for the development of a $q$-generalization of the well-known conjugate-gradient method, a research line currently being explored. \\

%%===============================================================================
\begin{acknowledgements}
The authors gratefully acknowledge the support provided by the National Counsel of Technological and Scientific Development (CNPq), Brazil.
\end{acknowledgements}
%%===============================================================================

%%===============================================================================

\end{document}